\documentclass[12pt]{article}
\usepackage[cp1251]{inputenc}
\usepackage[russian]{babel}
\usepackage{amssymb,amsthm,amsmath}

\theoremstyle{plain}

\newtheorem{Lemm}{Лемма}
\newtheorem{Theorem}{Теорема}
\newtheorem{Statement}{Утверждение}
\newtheorem{Note}{Замечание}

\newcommand{\suml}{\mathop{\sum}\limits}
\newcommand{\supl}{\mathop{\sup}\limits}
\renewcommand{\Re}{\operatorname{Re}}
\renewcommand{\Im}{\operatorname{Im}}

\renewcommand{\le}{\leqslant}
\renewcommand{\ge}{\geqslant}

\def\la{\lambda}
\def\eps{\varepsilon}

\def\C{{\mathbb C}}

\def\y{{\mathbf y}}

\def\A{{\mathbf A}}
\def\f{{\mathbf f}}
\def\xxi{{\mathbf \xi}}

\begin{document}

{\large {\bf \centerline{Равномерные оценки остатков}} {\bf
\centerline{в асимптотических выражениях для собственных функций}}
{\bf \centerline{оператора Штурма--Лиувилля с потенциалом --- распределением
\footnote{Работа поддержана грантом РФФИ No. 07-01-00283 и грантом ИНТАС No.
05-1000008-7883.}}} }

\bigskip

{\centerline{Савчук~А.~М.}}

\medskip

В настоящей заметке изучается оператор Штурма--Лиувилля
$$
L=-\dfrac{d^2}{dx^2}+q(x),
$$
в пространстве \(L_2[0,\pi]\) с граничными условиями Дирихле.
Предполагается, что потенциал \(q(x)\) является распределением
первого порядка сингулярности, т.е. \(q(x)\in W_2^{-1}\) или, что
то же самое, \(q(x)=u'(x),\) \(u(x)\in L_2[0,\pi]\) (производная
здесь понимается в смысле распределений). Операторы такого вида
были определены в работе \cite{SS1}. В работах \cite{SS1} и
\cite{SS2} было доказано, что оператор \(L\) фредгольмов с
индексами (0,0) (в случае вещественного потенциала --
самосопряжен), полуограничен, имеет чисто дискретный спектр. Были
получены асимптотические формулы для собственных значений и
собственных функций, а также была доказана базисность Рисса в
пространстве \(L_2\) системы собственных и присоединенных
векторов. В работе \cite{SS3} были получены асимптотические
формулы для собственных значений таких операторов, равномерные по
шару $\|u\|\le R$. В этой работе мы изучим поведение собственных и
присоединённых функций таких операторов, в частности получим для
них также равномерные оценки. Мы изучим случай потенциала,
принадлежащего шкале соболевских пространств
\(W_2^{-\theta}[0,\pi]\), где \(\frac12<\theta\le 1\).

Вначале мы дадим необходимые определения, введём обозначения и
приведём доказанные ранее (в работах \cite{SS1}, \cite{SS2} и
\cite{SS3}) результаты об операторе $L$.

Пусть в дифференциальном выражении
\begin{equation}\label{de1}
l(y)=-y''+q(x)y
\end{equation}
функция
 \(q(x)\in W_2^{-1}[0,\pi]\), а
\(u(x)=\int q(\xi)\,d\xi\) --- произвольная первообразная из
пространства \(L_2[0,\pi]\). Введём квазипроизводную
\[
    y^{[1]}(x)=y'(x)-u(x)y(x)
\]
и перепишем выражение~\eqref{de1} в виде
\begin{equation}\label{eq:1.1}
    l(y)=-\left(y^{[1]}\right)'-u(x)y^{[1]}-u^2(x)y.
\end{equation}
Несложно видеть, что для гладкой функции \(u(x)\) дифференциальное
выражение~\eqref{de1} и квазидифференциальное
выражение~\eqref{eq:1.1} совпадают. Теперь необходимо объяснить,
как мы понимаем решение уравнения
\begin{align*}
    l(y)&:=-y''+u'(x)y=\lambda y+f,&\lambda\in\mathbb C,\;f\in
    L_2,
\end{align*}
с начальными условиями $y(a)=y_0$, $y^{[1]}(a)=y_1$,
$a\in[0,\pi]$. Мы перепишем это уравнение в виде системы
\begin{equation*}
    \begin{pmatrix}y_1\\y_2\end{pmatrix}=\begin{pmatrix}u&1\\-\lambda-u^2
    &-u\end{pmatrix}\begin{pmatrix}y_1\\y_2\end{pmatrix}+\begin{pmatrix}
    0\\f\end{pmatrix},
\end{equation*}
где \(y_1=y\), \(y_2=y^{[1]}\). При этом элементы матрицы
\[
    \mathbf A(x)=\begin{pmatrix}u&1\\-\lambda-u^2&-u\end{pmatrix}
\]
являются функциями из \(L_1[0,\pi]\). Тогда (см., например,
\cite[гл. V, \S 16]{Na}) при любом $c\in[0,\pi]$ уравнение
\begin{align*}
\y'&=\A(x)\y+\f,& \y(c)&=\xxi\in\C^2,
\end{align*}
имеет единственное решение $\y(x)$, причем $\y(x)$ --- абсолютно
непрерывная на $[0,\pi]$ вектор--функция.

С выражением~\eqref{eq:1.1} свяжем оператор \(L\), определённый
равенствами
\begin{equation}\label{eq:L}
        Ly=l(y),\quad
        \mathfrak D(L)=\left\{y|\ y,y^{[1]}\in W_1^1[0,\pi],\,
        l(y)\in L_2[0,\pi],\,y(0)=y(\pi)=0\right\},
\end{equation}
Мы не предполагаем, что функция \(u(x)\) вещественна. Через
\(\overline{L}\) будем обозначать оператор, порождённый
сопряжённым дифференциальным выражением \(\overline{l}(y)\) (в
котором функция \(u(x)\) заменена на \(\overline{u(x)}\)) и теми
же краевыми условиями $y(0)=y(\pi)=0$. Прямым вычислением получаем
следующее предложение.

\begin{Statement}[Формула Лагранжа]\label{st:1} Для функций \(f\in\mathfrak
D(L)\), \(g\in\mathfrak D(\overline{L})\) справедливо тождество
\[
    (Lf,g)=(f,\overline{L}g),
\]
т.е. операторы \(L\) и \(\overline{L}\) взаимно сопряжены.
\end{Statement}

Представим основные результаты о спектре оператора $L$.
\begin{Statement}(см. \cite[Теорема 1.5]{SS2})\label{st:2}
Оператор \(L\) имеет непустое резольвентное множество и спектр его
дискретен.
\end{Statement}
Обозначим через $\omega(x,\la)$ решение дифференциального уравнения
$l(\omega)=\la\omega$ с начальными условиями $\omega(0,\la)=0$,
$\omega^{[1]}(0,\la)=1$. Ясно, что нули целой функции $\omega(\pi,\la)$
совпадают с собственными значениями оператора $L$. Алгебраической кратностью
собственного значения $\la_0$ мы будем называть кратность нуля $\la_0$ функции
$\omega(\pi,\la)$. С другой стороны, в силу теоремы существования и
единственности, геометрическая кратность каждого собственного значения равна
$1$. Пусть $\la_0$ есть собственное значение алгебраической кратности $p$.
Тогда функции $\omega^{(j)}_{\la}(x,\la)$, $j=1,2,\dots,p-1$ удовлетворяют
дифференциальным уравнениям
$l(\omega^{(j)}_{\la})=\la\omega^{(j)}_{\la}+\omega^{(j-1)}_\la$, причём
$\omega^{(j)}_{\la}(0,\la)=\omega^{(j)}_\la(\pi,\la_0)=0$. Тогда функции
$\omega^{(j)}_\la(x,\la_0)$, $j=0,1,\dots,p-1$ образуют цепочку собственной и
присоединённых функций, отвечающую собственному значению $\la_0$. Легко видеть,
что функции этой системы линейно независимы, а значит порождают подпространство
размерности $p$. Собственные значения мы обозначим $\{\la_n\}_1^\infty$, причём
нумерацию будем вести в порядке возрастания модуля $|\la_1|\le|\la_2|\le\dots$
с учётом алгебраической кратности. В случае совпадения модулей нумерацию будем
вести по возрастанию аргумента, значения которого здесь и далее договоримся
выбирать из полуинтервала $\arg\la_j\in(-\pi,\pi]$. Под системой собственных и
присоединённых функций $\{y_n(x)\}_1^\infty$ мы будем понимать систему,
полученную нормировкой $\|y_n\|_{L_2}=1$ системы
$\{\omega^{(j)}_\la(x,\la_n)\}$. Заметим, что для вещественной $u(x)$ все
собственные значения являются простыми (поскольку оператор $L$ здесь
самосопряжён и, стало быть, алгебраическая кратность собственных значений равна
их геометрической кратности).

Вернёмся к результатам о спектре оператора $L$. Символом
$l_2^\theta$ мы будем обозначать весовое пространство, состоящее
из последовательностей комплексных чисел $x=\{x_1,x_2,\dots\}$
таких, что
$$
\|x\|^2_\theta=\suml_{k=1}^\infty|x_k|^2k^{2\theta}<\infty.
$$
\begin{Statement}(см. Теорему 2.6 работы \cite{SS2})\label{st:3}
Пусть $u(x)\in W_2^\theta$ при некотором $0\le\theta<1/2$, а
$q(x)=u'(x)$ в смысле теории распределений. Тогда
\begin{equation}\label{evas}
\sqrt{\lambda_n}=n+s_n,\ n=1,2,\dots\,,\quad\text{где}\
\{s_n\}_{n=1}^\infty\in l^\theta_2.
\end{equation}
\end{Statement}
\begin{Note} В частности, отсюда следует, что числа $\sqrt{\la_n}$ лежат
в некоторой полуполосе $\{z\in\C\vert\,|\Im\, z|<\nu,\, \Re\, z>\mu>-\infty\}$
и все собственные значения являются простыми, начиная с некоторого номера $N$.
\end{Note}

В случае $\theta>0$ мы можем утверждать более сильный результат
\begin{Statement} (см. Теорему 2.1 работы \cite{SS3})\label{st:4}
Пусть $u(x)\in W_2^\theta$ при некотором $0<\theta<1/2$ и
$\|u\|_\theta\le R$, а $q(x)=u'(x)$ в смысле теории распределений.
Тогда
\begin{equation*}
\sqrt{\lambda_n}=n+s_n,\ n=1,2,\dots\,,\quad\text{где}\
\{s_n\}_{n=1}^\infty\in l^\theta_2,
\end{equation*}
причём $\|\{s_n\}\|_\theta\le C_{R,\theta}$, где величина $C$
зависит\footnote{Всюду в дальнейшем набор нижних индексов у постоянных,
встречающихся в различных оценках, будет обозначать множество тех и только тех
параметров, от которых зависит выбор данной постоянной.} лишь от $R$ и
$\theta$, но не зависит от функции $u$.
\end{Statement}
\begin{Note} Таким образом, при $\theta>0$ можно утверждать, что
числа $\sqrt{\la_n}$ лежат в некоторой полуполосе $\{z\in\C\vert\,|\Im\,
z|<\nu_{R,\theta},\, \Re\, z>\mu_{R,\theta}>-\infty\}$, а все собственные
значения являются простыми, начиная с некоторого номера $N_{R,\theta}$.
\end{Note}
Отметим, что в теореме 2.1 работы \cite{SS3} выписан и второй член
асимптотического разложения. В этой формуле также легко можно
получить равномерные оценки.

Перейдём к результатам о системе собственных и присоединённых
функций.
\begin{Statement} (см. теорему 2.9 работы \cite{SS2})\label{st:5}
Пусть $u(x)\in W_2^\theta$ при некотором $0\le\theta<1/2$. Тогда
система $\{y_n\}_1^\infty$ собственных и присоединённых функций
оператора \(L\) образует базис Рисса в пространстве $L_2[0,\pi]$.
\end{Statement}
Таким образом, существует единственная биортогональная система
$\{w_n\}_1^\infty$ (т.е. система, для которой выполнены равенства
$(y_n,w_m)=\delta_n^m$). В нашем случае эту систему можно
предъявить. Этот факт хорошо известен для классических потенциалов
$q\in L_1[0,\pi]$ (см., например, \cite[гл. 1, \S 3]{Ma}) и без
каких-либо изменений переносится на случай $q\in W_2^{-1}[0,\pi]$
(для доказательства требуется только формула Лагранжа (см.
утверждение \ref{st:1})).
\begin{Theorem}\label{th:1}
Пусть $u(x)\in L_2[0,\pi]$. Пусть $\la_n=\la_{n+1}=\dots=\la_{n+p_n}$ ---
собственное значение оператора $L$ кратности $p_n$. Обозначим
$\overline{\omega_{n+j}}=\overline{\omega}^{(j)}_\la(x,\la_n)$,
$j=0,1,\dots,p_n-1$ (нетрудно видеть, что это есть собственная и присоединённые
функции для сопряжённого оператора $\overline{L}$ и его собственного значения
$\overline{\la_n}$). Тогда функция $w_{n+j}$ биортогональной системы имеет вид
$$
w_{n+j}=\suml_{s=0}^{p_n-1}b_{j,s}\overline{\omega_s},
$$
где матрица $B=(b_{i,j})_0^{p_n-1}$ есть обратная матрица к
матрице $A=(a_{i,j})$, $a_{i,j}=(y_{n+i},\overline{\omega_j})$. В
частности, если $\la_n$ простое собственное значение, то
$w_n=\frac{\overline{y_n}}{(y_n,\overline{y_n})}$.
\end{Theorem}

Также, как и для собственных значений, асимптотическое поведение собственных
функций может быть легко определено. Асимптотически все собственные значения
однократны. Это позволяет найти асимптотическое поведение и векторов
биортогональной системы.
\begin{Theorem}\label{th:2}
Пусть $u(x)\in W_2^\theta[0,\pi]$, $0\le\theta<1/2$, а
$\|u\|_\theta\le R$. Тогда найдётся число $N_u$ такое, что
собственные функции $y_n$ и $w_n$ операторов $L$ и $\overline{L}$,
а также их квазипроизводные имеют асимптотику
\begin{equation}\label{efas}
\begin{array}{c}
y_n(x)=\sqrt{\frac2\pi}\sin nx+\phi_n(x),\qquad
y_n^{[1]}(x)=n\left(\sqrt{\frac2\pi}\cos nx+\phi^1_n(x)\right),\\
w_n(x)=\sqrt{\frac2\pi}\sin nx+\psi_n(x),\qquad
w_n^{[1]}(x)=n\left(\sqrt{\frac2\pi}\cos nx+\psi^1_n(x)\right),\ n=N,N+1,\dots.
\end{array}
\end{equation}
Остатки в этих формулах таковы, что последовательности чисел
$$\beta_n=\|\phi_n(x)\|_{C[0,\pi]}+\|\psi_n(x)\|_{C[0,\pi]}\quad\text{и}\quad
\gamma_n=\|\phi^1_n(x)\|_{C[0,\pi]}+\|\psi^1_n(x)\|_{C[0,\pi]}
$$
принадлежат $l_2^\theta$. При $\theta>0$ номер $N=N_{R,\theta}$, а
$\suml_{n=N}^\infty(\beta^2_n+\gamma^2_n)n^{2\theta}\le
C_{R,\theta}$.
\end{Theorem}
Этот результат мы приведём с доказательством, поскольку он не был
сформулирован в предыдущих работах. Нам потребуются некоторые
вспомогательные утверждения.
\begin{Lemm}\label{lem:1}
Пусть число \(\rho\in\mathbb C\) такое, что существует решение
\(\theta=\theta(x,\rho)\) уравнения
\begin{equation}\label{ekv:main}
    \theta(x,\rho)=\rho x+\int\limits_0^x u(t)\,
    \sin(2\theta(t,\rho))\,dt+\dfrac{1}{2\rho}\int\limits_0^x
    u^2(t)\,dt-\dfrac{1}{2\rho}\int\limits_0^x u^2(t)\,
    \cos(2\theta(t,\rho))\,dt.
\end{equation}
Тогда решение \(s(x,\rho)\) уравнения
$$
Ly=\rho^2y,
$$
удовлетворяющее условиям \(s(0,\rho)=0\), \(s^{[1]}(0,\rho)=1\),
допускает представление
\begin{equation}\label{polar}
    \rho s(x,\rho)=r(x,\rho)\,\sin\theta(x,\rho),\\
    s^{[1]}(x,\rho)=r(x,\rho)\,\cos\theta(x,\rho),
\end{equation}
где
\begin{equation}\label{r}
 r(x,\rho)=\exp\left\{-\int\limits_0^x u(t)\,
    \cos(2\theta(t,\rho))\,dt-\dfrac{1}{2\rho}\int\limits_0^x u^2(t)\,
    \sin(2\theta(t,\rho))\,dt\right\}.
\end{equation}
\end{Lemm}
\textit{Доказательство} этого утверждения имеется в п.~2.1
работы~\cite{SS2}. Функции $r$ и $\theta$ названы
модифицированными функциями Прюфера, а представление~\eqref{polar}
--- полярным представлением.
\begin{Lemm} (см. \cite[лемма 2.3]{SS3}\label{lem:2}
Пусть \(u(x)\in L_2\), \(\|u\|_{L_2}\leqslant R\), а \ $\nu$ "---
фиксированное число. Положим
\[
    \Upsilon(\rho)=\max\limits_{0\leqslant x\leqslant\pi}\left|
    \int\limits_0^x u(t)\sin(2\rho t)\,dt\right|+
    \left|\int\limits_0^xu(t)\,\cos(2\rho t)\,dt\right|+
    \dfrac{R^2(1+\varkappa+R\varkappa^2)}{2|\rho|},
\]
где \(\varkappa=\ch 2\pi\nu\). Тогда найдётся абсолютная
постоянная \(\varepsilon>0\) (можно взять \(\varepsilon=2^{-7}\)),
такая, что для всех \(\rho\), лежащих в полосе
\(|\Im\rho|\leqslant\nu\) и удовлетворяющих условию
\begin{align}\label{ekv:ineq}
    \Upsilon(\rho)&<\varepsilon(1+64R^2\varkappa^2)^{-2},
\end{align}
уравнение~\eqref{ekv:main} имеет единственное решение
\(\theta(x,\rho)\), которое представимо в виде
\begin{align}\label{ekv:pred}
    \theta(x,\rho)&=\rho x+f(x,\rho),&
    |f(x,\rho)|&\leqslant C\Upsilon(\rho),
\end{align}
где \(C\) --- абсолютная постоянная.
\end{Lemm}
\begin{Lemm}(см. \cite[лемма 2.4]{SS3}\label{lem:3}
Пусть \(u(x)\in W_2^{\theta}\), \(0\leqslant\theta\leqslant 1\),
\(\|u\|_{\theta}\leqslant R\). Тогда преобразование Фурье этой
функции
\begin{equation*}
    F(\rho)=\int\limits_0^{\pi} u(x) e^{i\rho x}\,dx
\end{equation*}
в полосе \(|\Im\rho|\leqslant\nu\) допускает оценку \(|F(\rho)|<CR
|\rho|^{-\theta}\), где \(C\) зависит только от \(\nu\).
\end{Lemm}
\begin{Lemm}\label{lem:4}
Пусть последовательность \(\{\rho_n\}_1^{\infty}\) такова, что
\(|\rho_n-n|<\delta<1/4\). Тогда при \(0\leq\theta<1/2\) оператор
\(T_x: W_2^{\theta}\to\ell_2^{\theta}\), определённый равенством
\begin{align*}
    T_xf&=\{c_n\}_1^{\infty},&c_n(x)&=\int\nolimits_0^x f(t)
    e^{i\rho_n t}\,dt,
\end{align*}
ограничен и его норма зависит только от \(\delta\) и \(\theta\).
\end{Lemm}
Итак, докажем теорему \ref{th:2}.
\begin{proof}
В силу леммы \ref{lem:1}
$\omega(x,\la)=s(x,\rho)=\frac{1}{\rho}r(x,\rho)\sin\theta(x,\rho)$.
Мы уже знаем, что все числа $\sqrt{\la_n}$ лежат в некоторой
горизонтальной полосе. Выберем число $\nu$ в лемме \ref{lem:2}
так, чтобы все они попали в полосу $\Pi_\nu=\{z\in\C\vert\,|\Im
z|\le\nu\}$. При этом $\nu=\nu(R,\theta)$, если $\theta>0$. Далее,
для любой функции $u$ имеем $\Upsilon(\rho)\to0$ при
$\Pi_\nu\ni\rho\to+\infty$, а значит условие \eqref{ekv:ineq}
выполнено при $\Re \rho>\rho_0(u)$. Более того, при $\theta>0$ мы
можем применить лемму \ref{lem:3}, а тогда $\Upsilon(\rho)\le
C_{R,\theta}|\rho|^{-\theta}$, т.е. условие \eqref{ekv:ineq}
выполнено при $\Re\rho>\rho_0(R,\theta)$. Из утверждения
\ref{st:4} следует, что вне полуполосы
$\{|\Im\rho|\le\nu,\,\Re\rho>\rho_0\}$ лежит лишь конечное число
$N=N_{R,\theta}$ собственных значений. Итак, приступим к
доказательству асимптотических равенств \eqref{efas}. Будем далее
считать что $\rho$ лежит в полуполосе
$\{|\Im\rho|\le\nu,\,\Re\rho>\rho_0\}$. Подставляя равенство
\eqref{ekv:pred} в тригонометрические функции, имеем
$$
\sin\theta(x,\rho)=\sin\rho x+\gamma_1(x,\rho),\
\sin2\theta(x,\rho)=\sin2\rho x+\gamma_2(x,\rho),\
\cos2\theta(x,\rho)=\cos2\rho x+\gamma_3(x,\rho),
$$
где $\supl_{x\in[0,\pi]}|\gamma_j(x,\rho)|\le
C_{\nu}\Upsilon(\rho)$. Здесь мы воспользовались следующим
известным фактом: \emph{если \(G(\xi)\) "--- аналитическая
функция, то \(|G(\xi)-G(\zeta)|\leqslant M\,|\xi-\zeta|\), где
\(M=\max|G'(\eta)|\) и максимум берётся по \(\eta\), лежащем на
отрезке \([\xi,\zeta]\) комплексной плоскости}. Подставляя
полученные равенства в представление \eqref{r}, получим
$r(x,\rho)=1+\gamma_4(x,\rho)$. Отсюда $s(x,\rho)=\frac{\sin\rho
x+\gamma_5(x,\rho)}\rho$, где
$\supl_{x\in[0,\pi]}|\gamma_j(x,\rho)|\le C_{\nu}\Upsilon(\rho)$.
Используем теперь асимптотические формулы \eqref{evas} для
собственных значений. Пусть вначале $\theta>0$. Для сокращения
записи будем обозначать через \(O_{R,\theta}\) последовательности
пространства $l_2^\theta$, норма которых оценивается величиной
$C_{R,\theta}$, зависящей от \(R\), \(\theta\), но не зависящей от
\(u\) в шаре \(\|u\|_{\theta}\leqslant R\). Поскольку
$\rho_n=n+s_n$, где $\{s_n\}=O_{R,\theta}$, то неравенство
$|s_n|<1/4$ выполнено при $n>N'_{R,\theta}$. Увеличивая, если
нужно, выбранный ранее номер $N$, и применяя лемму \ref{lem:4},
имеем $\Upsilon(\rho_n)=O_{R,\theta}$. Далее,
$$
\sin\rho_nx=\sin
nx+O_{R,\theta},\qquad\frac1{\rho_n}=\frac1n(1+O_{R,\theta})),
$$
а значит
$$
s(x,\rho_n)=\frac{\sin
nx}n+\frac1n\gamma_6(x,\rho_n),\quad\text{где}\
\supl_{x\in[0,\pi]}|\gamma(x,\rho_n)|=O_{R,\theta}.
$$
Теперь остаётся нормировать функции $s(x,\rho_n)$. Имеем
$\|s(x,\rho_n)\|_{L_2[0,\pi]}=\frac{\sqrt{\pi}}{n\sqrt{2}}+\frac1nO_{R,\theta}$,
откуда
$y_n(x)=\frac{s(x,\rho_n)}{\|s(x,\rho_n)\|_{L_2}}=\sqrt{\frac2\pi}\sin
nx+\phi_n(x)$, где
$\supl_{x\in[0,\pi]}|\gamma(x,\rho_n)|=O_{R,\theta}$. При
$\theta=0$ все выкладки повторяются с той разницей, что константы
теперь зависят от функции $u$.

Итак, первая асимптотическая формула в \eqref{efas} доказана.
Вторая формула доказывается абсолютно так же, поскольку
$s^{[1]}(x,\rho)=r(x,\rho)\cos\theta(x,\rho)$. Для доказательства
двух других равенств в \eqref{efas} заметим, что
$(y_n,\overline{y_n})=1+O_{R,\theta}$ и подставим это равенство в
формулу $w_n=\frac{\overline{y_n}}{(y_n,\overline{y_n})}$ (мы уже
выбрали выше $N$ таким, что $|\rho_n-n|<1/4$ при $n\ge N$, а
потому все наши собственные значения просты).
\end{proof}

Теперь мы сформулируем результаты о непрерывной зависимости
собственных значений и собственных функций от потенциала $u$.
Пусть $\la_n(u)$ --- простое собственное значение оператора
$L(u_0)$. Легко показать (см. \cite[лемма 5.3]{SS3}), что в
некоторой малой окрестности $\{u\in
L_2[0,\pi]\vert\,\|u-u_0\|<\eps\}$ оно является непрерывной и,
более того, аналитической функцией потенциала. То же верно и для
собственной функции $y_n(u)$, а из явного вида вектора
$w_n=\frac{\overline{y_n}}{(y_n,\overline{y_n})}$ биортогональной
системы вытекает непрерывная зависимость $w_n(u)$. Здесь нас,
однако, эти результаты не устраивают, поскольку
$\|w_n(u)\|\to\infty$, если функция $u$ изменяется так, чтобы
$|\la_n(u)-\la_m(u)|\to0$.
\begin{Theorem}\label{th:3}
Пусть $u\in W_2^\theta[0,\pi]$, где $0<\theta<1/2$, и
$\|u\|_\theta\le R$. Тогда существует такой номер $N_{R,\theta}$,
что для любого $n\ge N_{R,\theta}$ оператор проектирования
$P_n(u)$, определённый по правилу $P_nf=\suml_{k=1}^n(f,w_n)y_n$ и
действующий из пространства $L_2[0,\pi]$ в пространство
$W_2^1[0,\pi]$, непрерывно зависит от параметра $u$. А именно,
$\|P_n(u)-P_n(u_0)\|_{L_2\to W_2^1}\to0$, если
$\|u-u_0\|_\theta\to0$.
\end{Theorem}
\begin{Note}
Поскольку $\|P_n\|_{L_2\to L_2}\le\|P_n\|_{L_2\to W_2^1}$, отсюда
следует и непрерывная зависимость $P_n(u)$ как оператора из
пространства $L_2[0,\pi]$ в $L_2[0,\pi]$.
\end{Note}
Нам потребуется результат работы \cite{SS2} (теорема 1.9)
Обозначим через \(L_{\varepsilon}\) оператор, порождённый
дифференциальным выражением
\(l_{\varepsilon}(y)=-y''+q_{\varepsilon}(x)y\), где
$q_\eps=u_\eps'$.
\begin{Statement}\label{st:9}
Пусть $u_\epsilon\to u$ в $L_2[0,\pi]$. Существуют значения
\(\lambda\in\mathbb C\) такие, что при всех достаточно малых
\(\varepsilon\ge0\) значение \(\lambda\) принадлежит резольвентным
множествам операторов \(L_{\varepsilon}\), а последовательность
\((L_{\varepsilon}-\lambda)^{-1}\) сходится к $L$ при
\(\varepsilon\to 0\) в равномерной операторной топологии, т.~е.
$$
    \|(L_{\varepsilon}-\lambda)^{-1}-(L-\lambda)^{-1}\|\to 0
    \quad\text{при }\varepsilon\to 0.
$$
\end{Statement}
Теперь мы готовы доказать теорему \ref{th:3}.
\begin{proof}
В силу утверждения \ref{st:4} найдётся такое число
$\gamma_{R,\theta}=N_{R,\theta}+1/2$, что при любом $u$ из шара
$\|u\|_\theta\le R$ в области $\{\rho\in\C\vert\,|\Re
\rho|<\gamma_{R,\theta},\,|\Im \rho|<\gamma_{R,\theta}\}$ лежит ровно $N$ чисел
$\rho_k=\sqrt{\la_k}$, а для всех остальных чисел $\rho_k=\sqrt{\la_k}$
выполнены неравенства $|\rho_k-k|<1/4$. Для любого $n\ge N_{R,\theta}$
обозначим $\Gamma_n=\{\la\in\C\vert\,|\la|=(n+1/2)^2$ --- контур в
$\la$--плоскости. Тогда при любом $u$ из шара $\|u\|_\theta\le R$ внутри этого
контура лежит ровно $n$ собственных значений с номерами $k=1,2,\dots,n$. Теперь
наше утверждение следует из утверждения \ref{st:9} и классических результатов о
полунепрерывности изолированных частей спектра (см. \cite[теоремы IV.2.23 и
IV.3.16]{Ka}).
\end{proof}

\end{document}